\documentclass{amsart}
\usepackage{times,mathabx,amsfonts,amssymb,amsrefs,mathrsfs,tikz}
\usetikzlibrary{decorations,decorations.pathmorphing}

\usepackage[color]{showkeys}
\definecolor{gray}{gray}{0.5}

\newcommand\Z{\mathbb{Z}}

\newtheorem{lemma}{Lemma}[section]

\newtheorem{proposition}[lemma]{Proposition}
\newtheorem{theorem}[lemma]{Theorem}
\newtheorem{corollary}[lemma]{Corollary}
\newtheorem{example}[lemma]{Example}

\newtheorem{conjecture}[lemma]{Conjecture}

\theoremstyle{definition}
\newtheorem{remark}[lemma]{Remark}
\newtheorem{definition}[lemma]{Definition}
\newtheorem{question}[lemma]{Question}

\font\manfnt=manfnt

\newcommand\twoheaddownarrow{\hbox to 0pt{\raisebox{0.3ex}{$\downarrow$}}
  \hbox to 0pt{\raisebox{-0.2ex}{$\downarrow$}}\phantom\downarrow}

\begin{document}
\title{Almost invariance of  distributions for random walks on groups}

\author{Anna Erschler}
\address{CNRS, DMA, \'{E}cole Normale Sup\'{e}rieure, 45 rue d'Ulm, Paris, France}

\begin{abstract}
We  study  the neighborhoods of a typical  point  $Z_n$ visited at $n$-th step of a random walk, determined by the condition that the transition probabilities stay close to $\mu^{*n}(Z_n)$.  
If such neighborhood contains a ball of radius $C \sqrt{n}$, we say that the random walk has almost invariant transition probabilities.
We prove  that simple random walks on  wreath products of $\Z$ with  finite groups have almost invariant distributions. A weaker version of almost invariance implies a necessary  condition of Ozawa's criterion for the property $H_{\rm FD}$.
We define and study the radius of almost invariance, 
  we estimate this  radius for random walks on iterated wreath products and show this radius can be asymptotically strictly smaller than $n/L(n)$, where $L(n)$ denotes the drift function of the random walk. We show that the radius of individual almost invariance of a simple random walk on the wreath product of $\Z^2$ with a finite group is asymptotically strictly larger than $n/L(n)$. Finally, we show the existence of groups such that the radius of almost invariance is smaller than a given function, but remains unbounded. We also discuss possible limiting distribution of ratios of transition probabilities on non almost invariant scales.

\end{abstract}

\maketitle

\section{Introduction}

Given a sequence of functions $f_n:G \to \mathbb{R}$, a sequence of probability measures $\mu_n$, and
a sequence of subsets $\Omega_n \subset G$ , we can consider  the
ratios
  $f_n(xg_n )/ f_n(x)$, where $g_n\in \Omega_n$) and ask whether this ratio
is close to one for $x$ chosen with  $\mu_n$ probability close to one.
  If the latter happens for $\Omega_n$ being the ball of of radius  $r(n)$, centered at the identity, we say that $f_n$ is
 almost invariant on the scale $r(n)$.

Otherwise, if this is not the case, we can ask whether $f_n( x g_n) / f_n(x)$   is close, with probability one with respect to $\mu_n$, to a given probability distribution?

In this paper, we study the case when both $f_n$ and $\mu_n$ are $n$ step transition probabilities of a random walk $(G,\mu)$. Here $G$ is some finitely generated group and $\mu$ is a probability measure on $\mu$, $f_n(g) =\mu_n(g)=\mu^{*n}(g)$.

\begin{definition} \label{almostinvariantscales}
Let $G$ be a group, generated by a finite set $S$.
Random walk $(G,\mu)$ has {\it almost invariant distributions} on the scale $r(n)$, if for all $\epsilon>0$ and all $n$ there exists a subset $V_{\epsilon,n} \subset G$, satisfying  $\mu^{*n}(V_n) \ge 1-\epsilon$, such that the following holds.
For any sequences $h_n \in V_n$  and $g_n\in G$, such that $l_{G,S}(g_n) \le r(n$) it holds $|\mu^{*n}( h_n g_n)/\mu^{*n}(h_n)-1| \le \epsilon$.
\end{definition}

For example, non-degenerate aperiodic random walks on  finite groups are almost invariant on any scale.
The estimates of Hebisch and  Saloff-Coste imply that  simple aperiodic random walks on nilpotent groups are almost invariant on any  scale smaller than $\sqrt{n}$ (see Example \ref{examplenilpotent}).

\begin{definition} \label{definitionradiusalmostinvariance}
Given a function $r(n)$, we say that {\it  the radius of almost invariance} for transition probabilities of the random walk $(G,\mu)$ is asymptotically larger than $r(n)$ if $\mu^{*n}$ is almost invariant on the scale  $r''(n)$ for any $r''(n)$ such that $r''(n)/r(n) \to 0$.  We write in this case $r_{\rm a.i.}(n) \succeq r(n)$.
We say that the radius of almost invariance is asymptotically smaller than $r(n)$ if for any function $r'(n)$ such that $r'(n)/r(n)$ tends to infinity on some subsequence,
   $\mu^{*n}$ is not almost invariant on the scale $r'(n)$.  We write in this case  $r_{\rm a.i.}(n) \preceq r(n)$.
If $r_{\rm a.i.}(n) \succeq r(n)$ and $r_{\rm a.i.}(n) \preceq r(n)$, we say that $r_{\rm a.i.}(n) = r(n)$ is the radius of almost invariance for transition probabilities of the random walk $(G,\mu)$.

\end{definition}

We believe that  the radius of almost invariance is at most $\sqrt{n}$  for any simple random walk with the support generating an infinite group.
 We say that the random walk $(G,\mu)$ {\it has almost invariant transition probabilities}  if it is almost invariant on the  hypothetically largest possible scale: that is, 
 if $r_{\rm a.i.}(n) \succeq \sqrt{n}$.

We would like to stress that it is essential for the questions we study in this paper that the centers of neighborhood are chosen at random with respect to $\mu^{* n}$. Ratios of transition probabilities in the neighborhoods of a deterministic point have been studied in various situations: in the amenable case where one can study what we call {\it radii of almost constancy} (see the definition below and section \ref{fixedpoint}); in the non-amenable case where one can study {\it ratio theorems}  for transition probabilities. Ratio $\mu^{*n}(hg)/\mu^{*n}(h)$, with $h$ chosen at random on the trajectory of the random walk $(G,\mu)$ have been studied so far only for fixed increments, according to our knowledge.  In the latter case, it has been proven by Kaimanovich and Vershik \cite{kaimanovichvershik} that closeness of such ratios to one characterizes non-triviality of  the Poisson boundary.

The radius of almost invariance describes the size of a neighborhood of a typical point (chosen at random with respect to $\mu^{*n}$) such that the transition probabilities are almost constant in such neighborhood. If instead of
 a random point we consider a fixed point, the size of such neighborhood can be quite different, 
see Subsection \ref{fixedpoint}). 
The existence of a non-trivial neighborhood of the identity with almost constant decay for a symmetric random  walk on $G$  is equivalent to the amenability of $G$ (as it follows from the ratio theorem of Avez, see  Lemma \ref{lemmaamenable}), while the existence of non-trivial scale of almost invariance is equivalent to the triviality of the Poisson boundary (as it follows from the results of Kaimanovich and Vershik, mentioned above, see Proposition \ref{propositiontrivialboundary}). 

Let us say that $R(n)$ is the {\it radius of almost constancy} for transition probabilities of the random walk $(G,\mu)$ if there exists $C>0$ such that $\mu^{*n}(g)/\mu^{*n}(e)$ is bounded (by $C$ from above and $1/C$ from below) for any element $g$ of length at most $R(n)$ and if for any increasing sequence $n_i$ and $R'$  such that $R'(n_i)/R(n_i)$ tends to infinity, there exists a sequence $g_i$, of length at most $R'(n_i)$, such that $\mu^{*{n_i}}(g)/\mu^{*{n_i}}(e)$ is not bounded. This definition is mostly interesting for symmetric measure $\mu$, and in this case  the ratio in question is bounded if it is bounded from below, for even $n$.

The radii of almost invariance  and almost constancy of transition probabilities in the neighborhood of a fixed point are asymptotically equivalent for simple random walks on nilpotent groups.  For all random walks with non-trivial boundary
(as well as for various examples of random walks with trivial boundary) the radius of almost constancy in the neighborhood of a fixed point is  asymptotically equal or larger than the radius of almost invariance.
We show however that  there exist groups such that the radius of almost invariance for  simple random walks on such groups is larger than the radius of almost constancy in the neighborhood of a fixed point. This is the case for wreath product of $\Z$ with a finite group, as follows from Theorem \ref{theorem1}.

Recall that the {\it wreath product} of $A$ and $B$, which we denote by $A\wr B$,  is a semi-direct product
of $A$ and $\sum_A B$, where $A$ acts on $\sum_A B$ by shifts.

A probability measure $\mu$ on  the wreath product of $A$ with $B$,   defines {\it a switch-walk-switch random walk}, if $\mu=\mu_B*\mu_A*\mu_B$, where $\mu_A$ and $\mu_B$ are measures on $A$ and $B$ correspondingly. 

In case when $A=\mathbb{Z}$ and $B$ is a finite group, we say that $\mu$ defines a standard simple switch-walk-switch random walk on $G = A \wr B$,  if $\mu=\mu_B*\mu_A*\mu_B$, where $\mu_A$ is a symmetric measure on $\mathbb{A }$ with its support being equal to ${-1,0,1}$ and
$\mu_B$ is the equidistribution on $B$.

\begin{theorem}  \label{theorem1}

 Let $G= \mathbb{Z} \wr F$, $F$ is a finite group. Let $(G, \mu)$ be a standard simple switch-walk-switch random walk on $G$.

{\rm 1) [Almost Invariance.] } For any $\epsilon$ there exists $c>0$  
and asubsets 
$V_n$, such that $\mu^{*n}(V_n)\ge 1-\epsilon$ and such that for any $h \in V_n$ and any
$g \in B(e, c\sqrt{n})$
$$
|\mu^{*n}(hg)/\mu^{*n}(h)-1| \le \epsilon,
$$
where $B(e, cn)$ is the ball of radius $cn$ is some  word metric of $G$.

{\rm 2) [Invariance]} For any $\epsilon$ there exists $c>0$  
and subsets 
$V_n$, such that $\mu^{*n}(V_n)\ge 1-\epsilon$ and such that for any $h \in V_n$ and any
$g \in B(e, c\sqrt{n}) \cap G_0$
$$
\mu^{*n}(hg)=\mu^{*n}(h),
$$
where $G_0= \sum F \subset G$.

\end{theorem}

The first claim of Theorem \ref{theorem1} shows that the radius of almost invariance is $\succeq \sqrt{n}$, and this implies that 
 $G= \mathbb{Z} \wr F$ satisfies the first assumption of proposition in Section 4 of \cite{ozawa}. 
(In fact, a weaker claim that the {\it radius of individual almost invariance} (see Definition \ref{definitionradiusalmostinvariance}
 in Section \ref{preliminaries}) is asymptotically larger than $\sqrt{n}$ is sufficient for this conclusion).
 Since it  is well known and easy to see that the wreath products of $\Z$ with finite groups admit a sequence of subsets $U_n \subset B(e,CN)$ for some $C>0$, such that $\# \partial U_n/\# U_n \le 1/n$ (that is, they admit {\it a controlled Foelner sequence} in the terminology of \cite{cornuliertesseravalette}), these wreath products satisfy also the second assumption of Ozawa's proposition. 
Theorem \ref{theorem1} provides therefore an answers for  a question of Ozawa (private communication and a remark before the proposition in Section 4 of \cite{ozawa}) about existence of groups of super-polynomial growth where his criterion can be applied to prove the property $H_{FD}$, defined by Shalom in \cite{shalom}. We remind that the property $H_{FD}$ states that any unitary representation with non trivial first reduced cohomology group 
of a group $G$ on a Hilbert space contains a finitely dimensional invariant subspace. Since any group without property T of Kazhdan , in particular any amenable group, admits a unitary representation with non-trivial first reduced cohomology group  (\cite{mok, korevaarschoen}, see also \cite{shalom} and appendix in \cite{ozawa}), this property provides a sufficient criterion for an amenable group to admit a non-trivial virtual quotient to $\Z$.

Construction of virtual homomorphism to $\Z$ is the key step in all known proofs of the Polynomial Growth theorem of \cite{gromov}.
We remind that the original argument by Gromov \cite{gromov} uses the fact that the space, obtained as a limit of normalized word metrics of $G$, is locally compact, and this is never the case 
when the growth is super-polynomial. The proof of Kleiner \cite{kleiner} uses the fact that the space of $\mu$-harmonic functions of linear growth  on a group of polynomial growth  is finitely dimensional ,
and a natural conjecture is that the space of  $\mu$-harmonic functions of linear growth  on any group of super-polynomial growth
is infinite dimensional. In contrast with these proofs, we show that the argument of Ozawa, obtained in his proof of Polynomial growth theorem, can be applied to some groups of super-polynomial growth.

While the assumption on the measure in  Theorem \ref{theorem1} is important for the second claim, this assumption is not necessary for the first claim of the theorem. It can be shown that the same conclusion holds for any simple random walk.
Moreover, the argument in the proof of the first part of Theorem \ref{theorem1} can be applied to large classes of  groups, we plan to return to this question elsewhere. We ask in this context:

\begin{question} Let $G$ be a group of finite Pruefer rank (e.g. polycyclic) and $\mu$ be a non-degenerate symmetric finitely supported measure on $G$. Is it true that the distributions of $(G,\mu)$ are almost invariant?
\end{question}

Recall that the drift function $L(n)$ is the mean displacement after $n$ steps of the random walk, that is,
$L(n)= \sum_{g\in G} l_S(g) \mu^{*n}(g)$.

In contrast with $L(n)$, the radii of almost invariance are defined for measures not necessarily satisfying a moment condition. The asymptotic behavior of the radii of almost invariance is related to, but not defined by the drift function. 
 
 In Proposition \ref{propositioniteratedwreath} we show that the almost invariance radius of a simple switch-walk-switch random walk on 
$j$ times iterated wreath product of $\Z$ satisfies

 $$
 r_{\rm a.i.}(n) \preceq C   n^{1/2^{i+1}}/ \sqrt{\ln(n)}.
 $$
  This provides examples of simple random walks with almost invariance radius strictly smaller than 
 $n/L(n)$, since in this example $n/L(n) \sim  n^{1/2^{i+1}}$. 
 
In Proposition \ref{propositionwrz2}
we show that the radius of individual almost invariance (as defined in definition Definition \ref{definitionradiusalmostinvariance}) of a standard switch-walk-switch random walk on a wreath product of $\Z^2$ with a finite group is bounded from below by a radius of a ball covered by its projection to $\Z^2$. This implies that the radius of individual almost invariance satisfies for all sufficiently large $n$
$$
 \bar{r}(n) \succeq \exp(\sqrt{\ln(n)})
$$
 (see Corollary \ref{corollary2lowerboundz2}), and thus  it is asymptotically larger than $n/L(n)$, since in these examples $n/L(n) \sim \ln(n)$.

\begin{question}
Is it true that any group of finite Pruefer rank satisfies  the property $H_{FD}$ of Shalom? Does any group such that the drift of a simple random walk on this group is asymptotically equivalent to $\sqrt{n}$ satisfy the property $H_{FD}$ of Shalom?

\end{question}

Since groups with property $H_{FD}$  admit a finite index subgroup with infinite Abelinisation (\cite{shalom}),
the positive answer to the question above would imply two conjectures below.

\begin{conjecture}
The drift function of any simple random walk on an infinite finitely generated simple group is strictly larger than $\sqrt{n}$.
\end{conjecture}

\begin{conjecture}
The drift function of any simple random walk  on an infinite finitely generated torsion group is strictly larger than $\sqrt{n}$.
\end{conjecture}

It is proven in \cite{erschlercritical}, Corollary 1 that the drift of simple random walks on the first Grigorchuk group satisfies $L(n) \ge n^\gamma$, for some $\gamma >1/2$ and infinitely many $n$.
 It can be shown that the  drift of torsion groups, e.g. some of Grigorchuk groups from \cite{grigorchukdegrees} that are close on some scales to  solvable groups,
can be arbitrary  close to $\sqrt{n}$ on an infinite subsequence (moreover, on a sequence of  arbitrarily quickly growing intervals).

It seems interesting to describe smallest possible  drifts for simple random walks on infinite simple groups and on infinite torsion groups.


The paper has the following structure. In Section \ref{preliminaries} we describe first examples and basic properties of almost invariance and radii of almost invariance. In Proposition  \ref{propositiontrivialboundary}  we describe the relation with the Poisson boundary. We define a weaker notion of {\it individual almost invariance} on a given scale, and 
 in Corollary  \ref{corollaryentropy} we give a lower bound for the  radii of individual almost invariance in terms of entropy of random walks.
    Lemma \ref{lemmaconditionned} is essential for optimal lower bounds for radii of individual almost invariance for various groups: the idea is that in order to prove almost invariance for $\mu^{*n}$, it is sufficient to chose a sequence of events $\alpha_n$ on the space of the trajectories of the random walk $(G,\mu)$, such that the probabilities of $\alpha_n$ are close to one, and to prove the almost invariance for the sequence of conditional measures $\mu_n=\mu^{*n}(g)| \alpha_n$.     A similar ideas are used in the proof of Theorem \ref{theorem1} and Proposition \ref{propositionwrz2} to get a stronger conclusion and to bound from below (non-individual) radius of almost invariance.
In Subsection \ref{fixedpoint} we compare the radius of almost invariance with that of almost constancy of transition probabilities.

In Section \ref{prooftheorem1} we prove Theorem \ref{theorem1}.

In Section \ref{morewreath} we estimate the radii of almost invariance for distributions of simple random walks on iterated wreath products with $\mathbb{Z}$ and wreath products of $\mathbb{Z}^2$ with finite groups.

In Section \ref{piecewise} we prove the following proposition that shows that invariance radii 
can be arbitrary small.

\begin{proposition} \label{propositionbadinvariance}
For any $F(n)$ tending to infinity there exists a finitely generated group $G$ such that the almost invariance radius of any simple random walk on $G$ is not bounded 
 and $r_{\rm a.i.}(n) \preceq F(n)$.
 

\end{proposition}

In the last Section we show that,
though no group with non-trivial boundary can admit a non-trivial scale for almost invariance, it may happen
that $\mu^{*n}(hg_n)/\mu^{*n}(h)$, normalized by a constant $C_g$, admits a (non-constant) limit distribution,
see Example \ref{examplefree} for the case of standard simple random walks on a free non-Abelian group.

{\bf Acknowledgements. } I am grateful to Vadim Kaimanovich for many helpful discussions and comments on the preliminary version of this paper. I am also grateful to Balint Toth for turning my attention to
the result of \cite{demboperesrosen}.

\section{Preliminaries} \label{preliminaries}

Recall that the period of  a probability measure $\mu$ on  a group $G$ is the greatest common divisors of $K\ge0$ such that $\mu^{*K}(e)>0$. A probability measure $\mu$ is called {\it aperiodic} if its period is equal to one. For example, if $\mu(e)>0$ then $\mu$ is aperiodic.
Observe that for any symmetric measure on a group it holds $\mu^{*2}(e)>0$, and hence 
a symmetric measure $\mu$ on a group $G$ is aperiodic if and only if there exists an odd integer $K$ such that $\mu^{*K}(e) \ne 0$.

\begin{example} [Finite groups] \label{examplefinite}
Let $F$ be a finite group,  and $\mu$ be a non-degenerate aperiodic random walk on $F$. For any sequence $r(n)$ the distributions of $(G,\mu)$ are almost invariant
on the scale $r(n)$.
\end{example}

\begin{example} [Nilpotent groups]\label{examplenilpotent}
Let $G$ be a  finitely generated virtually nilpotent group,  $\mu$ be a finitely supported  symmetric non-degenerate aperiodic random walk on $F$.
 Then the distributions of the random walk $(G,\mu)$ are almost invariant on the scale  $r(n)$ if $r(n)/\sqrt{n} \to 0$. 
\end{example}

{\bf Proof.} Follows from  the "gradient estimate" of Hebisch and Saloff-Coste (see line 2 on p.675 and (14) of Theorem 5.1 in  \cite{hebischsaloff}), who prove, under assumption $\mu(e)>0$  that
$$
\mu^{*n} (g) -\mu^{*n}(gs) \le C' n^{-(D+1)/2} \exp(-l_S ^2(x)/CN),
$$
where  $S$ is a finite generating set of $G$ (see first line on page 675, as well as (14) of Theorem 5.1 ).

This estimate shows that if $l_S(g) \le K \sqrt{n}$, then $\mu^{*n} (g) -\mu^{*n}(gs) \le C'' n^{-(D+1)/2} $ for some $C''>0$ and any $s \in S$.
Therefore,   $\mu^{*n} (g) -\mu^{*n}(gh) \le C''n^{-D/2} l(h)/\sqrt{n}   \le  C''K \epsilon n^{-D/2}   $ for any $g: l_S(g) \le K \sqrt{n}$ and $h:l_S(h) \le \epsilon \sqrt{n}$.   Since
$\mu^{*n}(g) \ge C_0  n^{-D/2}$  for some $C_0>0$ depending on $K$ and any  $g: l_S(g) \le K \sqrt{n}$ (\cite{hebischsaloff}), we conclude that $|\mu^{*n}(g)/\mu^{*n}(gh) -1|  \le \epsilon K''  $ for some $K''$ depending on $K$, and all
$g: l_S(g) \le K \sqrt{n}$ and $h:l_S(h) \le \epsilon \sqrt{n}$.

Therefore, for each $K$ there exists $\delta>0$ such that the following holds. Let  $r(n)$ satisfies $r(n)/\sqrt{n} \to 0$. Consider the ball$V_n= B_{G,S}(e,K \sqrt{n})$. Then for any $h_n \in V_n$ and any $g_n: l_{G,S}(g_n)\le r(n)$ it holds
$$
|\mu^{*n}(h_ng_n)/\mu^{*n}(h_n)-1| \le \delta.
$$


Since the statement holds true for measures with $\mu(e)\ne e$, it is clear it holds also for any aperiodic measure.


 Definition \ref{almostinvariantscales} of almost invariance on some scale requires the existence of one single set $V_n$ for all increments of length at most $r(n)$. In the Definition below we give a weaker version of this n notion, where $V_n$ can depend on the sequence $g_n$ satisfying $l_{G,S}(g_n) \le r(n)$.

\begin{definition} \label{definitionalmostinvariantscaleindividually}
Let $G$ be a group, generated by a finite set $S$.
Random walk $(G,\mu)$ has {\it individually almost invariant distributions} on the scale $r(n)$, if for all $\epsilon>0$  and all $g_n$ satisfying $l_{G,S}(g_n) \le r(n)$
 there exist subsets $V_{\epsilon,n} \subset G$, such that   $\mu^{*n}(V_n) \ge 1-\epsilon$ and such that 
for any sequences $h_n \in V_n$  it holds $|\mu^{*n}( h_n g_n)/\mu^{*n}(h_n)-1| \le \epsilon$.
\end{definition}

If in the definition above we replace $|\mu^{*n}( h_n g_n)/\mu^{*n}(h_n)-1| \le \epsilon$ by   the condition that $\mu^{*n}( h_n g_n)/\mu^{*n}(h_n)-1$ is bounded from below and from above by some positive constant and the condition $\mu^{*n}(V_n) \ge 1-\epsilon$ by $\mu^{*n}(V_n) >a$ for some positive constant $a$ and all $n$, then we say that the random walk has individually almost invariant distributions weakly on the scale $r(n)$.

\begin{definition} \label{definitionradiuindividualsalmostinvariance}

Given a function $r(n)$, we say that {\it  the radius of individual almost invariance} for transition probabilities of the random walk $(G,\mu)$ is asymptotically larger than $r(n)$ if $\mu^{*n}$ is individually almost invariant on the scale  $r''(n)$ for any $r''(n)$ such that $r''(n)/r(n) \to 0$.  We write in this case $\bar{r}(n) \succeq r(n)$.
We say that the radius of iindividual almost invariance is asymptotically smaller than $r(n)$ if for any function $r'(n)$ such that $r'(n)/r(n)$ tends to infinity on some subsequence,
   $\mu^{*n}$ is not almost invariant on the scale $r'(n)$.  We write in this case  $\bar{r}(n) \preceq r(n)$.
If $\bar{r}(n) \succeq r(n)$ and $\bar{r}(n) \preceq r(n)$, we say that $\bar{r}(n) = r(n)$ is the radius of almost invariance for transition probabilities of the random walk $(G,\mu)$.

\end{definition}

\begin{remark} [General Markov chains] For random walks on groups, it is not important whether we consider multiplication on the left or from the write in the definition of almost invariance. For general random walks and individual sequence $g_n$, the definition makes sense for the multiplication on the left and there
seems to be no analog of this notion in the case of multiplication on the right. However, for uniform almost invariance on a given scale, one can define this notion for the multiplication on the right as well.

\end{remark}

\begin{proposition} [Relation with the triviality of the Poisson boundary] \label{propositiontrivialboundary}

Let $G$ be a finitely generated group and $\mu$ be a non-degenerate aperiodic measure on $G$.  The following properties are equivalent

\begin{enumerate}

\item
There exists a   function
$R(n)$,  tending to infinity, such that the distributions of the random walk $(G,\mu)$ are s almost invariant on the scale $R(n)$.


\item
There exists a sequence
$R(n)$,  tending to infinity on some subsequence, such that the distributions of the random walk  $(G,\mu)$ are weakly individually almost invariant distributions on the scale $R(n)$.

\item
The Poisson boundary of $(G,\mu)$ is trivial.

\end{enumerate}
\end{proposition}

{\bf Proof.} It is clear that $(1)$ implies $(2)$.

Let us show that $(2)$ implies $(3)$. Take a random walk $(G,\mu)$ with non trivial Poisson boundary.

By the Law of $0$ and $2$ (see Kaimanovich,  \cite{kaimanovich02law}, Theorem 2.1, p.155), we know that for any $\epsilon>0$  there exists probability distributions such that $\sum _{g\in G}|\nu_1 \mu^{*n}(g) -\nu_2 \mu^{*n}(g)| \ge 2-\epsilon$ for all sufficiently large $n$. Observe that for any $\epsilon>0$ we can chose $\nu_1$ and $\nu_2$ as above  to have finite support. Observe also that if  $\sum _{g\in G}|\nu_1 \mu^{*n}(g) -\nu_2 \mu^{*n}(g)| \ge 2-\epsilon$ for some $n$, then for some $h_1$ and $h_2$ in the supports of $\nu_1$ and $\nu_2$ it holds
 $\sum _{g\in G}|h_1 \mu^{*n}(g) -h_2 \mu^{*n}(g)| =\sum _{g\in G}| \mu^{*n}(g) -h^{-1}h_2 \mu^{*n}(g)| 
 \ge 2-\epsilon$. Since $\sum _{g\in G}| \mu^{*n}(g) -h^{-1}h_2 \mu^{*n}(g)| $ is a non-increasing function in $n$,
 this implies that for any $\epsilon >0$ there
exists $h\in G$ such that 
 $\sum _{g\in G}| \mu^{*n}(g) -h \mu^{*n}(g)| \ge 2-\epsilon$ for all sufficiently large $n$.

 Therefore, for any $\epsilon>0$ and $C>0$  there exists $h\in G$ such that  for all sufficiently large $n$.
$$
\mu^{*n}(g \in G: \mu^{*n}(g h)/\mu^{*n}(g) >C)\ge 1-\epsilon.
$$
This is in contradiction with $(2)$.

Now we show that $(3)$ implies $(1)$.  
Since the Poisson boundary is trivial, we know from Theorem 4.2 of  \cite{kaimanovichvershik} that for any $g\in G$ 
$$
\mu^{*n}(x \in G:  |1-\mu^{*n}(g x)/\mu^{*n}(x)|>\epsilon ) \to 0,
$$
as $n$ tends to $\infty$.

Take $A: \exp(An) \ge v_S(n)$ for all $n$, where $v_S(n)$ is the growth function of $G$ with respect to $S$.
Observe that for any  $N_1<N_2 <N_3 \dots... \in \mathbb{N}$  there exist $M_1<M_2<\dots$ such that for any $n>M_i$ and any $g: l_S(g) \le N_i$
$$
\mu^{*n}(x\in G:  1- 1/N_i <\mu^{*n}( xg)/\mu^{*n}(x) < 1+1/N_i ) > 1-\exp(-A N_i)/N_i.
$$

Note that there exists $V_i$ such that for all $n\ge M_i$ it holds $\mu^{*n}(V_i) \ge 1/N_i$ and 
$$ 
1- 1/N_i <\mu^{*n}(x g)/\mu^{*n}(x) < 1+1/N_i
$$
for all $g: l_S(g) \le N_i$ and all $n\ge M_i$.

Take a non-decreasing function  $r(n)$, tending to $\infty$ as $r$ tends to $\infty$  such that $r(M_i) \le {N_{i-1}}$ for all $i$. It is clear that  
for any $\epsilon$, any sufficiently large $n$ and any $g$ of length at most $r(n)$ 
$$
\mu^{*n}(x \in G:  |1-\mu^{*n}(g x)/\mu^{*n}(x)|>\epsilon ) \le \epsilon.
$$

This completes the proof of the proposition.

\begin{remark} \label{remarkreiter} Let $\mu$ and $\mu'$ be probability measures on a space $X$.

1) If for some $K, \epsilon>0$ $\mu(g)/\mu'(g) \le K$ with probability $\ge 1-\epsilon$  with respect to
$\mu$, then 
$$
\sum_{x\in X} |\mu(x) - \mu'(x)| \le 2( \epsilon +K-1)
$$

2) If for some $K_2, \epsilon_2>0$ it holds $\mu(g)/\mu'(g) \ge K_2$ with probability $\ge \epsilon_2$  with respect to
$\mu$, then
$$
\sum_{x\in X} |\mu(x)-\mu'(x)| \ge \epsilon_2 (K-1)
$$

{\bf Proof.}

1).  Observe that 

$$
\sum_{x\in X} |\mu(x) - \mu'(x)|  = 2  \sum_{x\in X:  \mu(x)\ge \mu'(x)   } |\mu(x) - \mu'(x)|          \le 2( \epsilon +K-1)
$$

2). The proof is straightforward.

\end{remark}

\begin{corollary} The transition probabilities of the random walk $(G,\mu)$ are  almost invariant with respect to a sequence $g_n$ if and only if the total variance between convolutions of $\mu$ and shifted convolutions of $\mu$ tends to zero: 
$$
\sum_{g\in G} |\mu^{*n}(g) -  \mu^{*n}(g g_n)| \to 0,
$$
as $n \to \infty$.
\end{corollary}

\begin{remark} \label{remarkdifmeasures}Let $S$ be a finite subset of $G$. Suppose that for any $s\in S$ it holds
$$
\sum_{h}|(\mu^{*n}(h)   -\mu^{*n}(hs)| \le 1/F(n).
$$
Then for any $g$ in the subgroup generated by $S$ 
$$
\sum_{h}|(\mu^{*n}(h)   -\mu^{*n}(hg)| \le l_S (g) /F(n).
$$

\end{remark}

{\bf Proof.}
Observe that
$$
\sum_{h}|(\mu^{*n}(h)   -\mu^{*n}(hab)| \le \sum_{h}|(\mu^{*n}(h)   -\mu^{*n}(ha))| + \sum_{h}|(\mu^{*n}(h)   -\mu^{*n}(hb))|.
$$

\begin{corollary} \label{corollary1} Let $\mu$ be a finitely supported measure on  $G$  such that 
for some $F(n)$ and 
for any $s$ in the support of $\mu$
$$
\sum_{h}|(\mu^{*n}(h)   -\mu^{*n}(hs)| \le 1/F(n).
$$

Take $f(n)$ such that $f(n)/F(n) \to 0$ as $n\to \infty$. 
The distributions of $(G,\mu)$ are individually almost invariant on the scale $f(n)$.

\end{corollary}
 
{\bf Proof.} Follows from  Remark \ref{remarkdifmeasures} and 2) of  Remark \ref{remarkreiter}.

Recall that an {\it entropy} of a probability measure $\nu$ on a space $X$ is $-\sum_{x}\in X  \nu(x)\log(\nu(x))$. Given a probability measure $\mu$ on a group $G$ denote by $H(n)$ the entropy of its $n$-th convolution, $H(n)= H(\mu^{*n})$.

\begin{corollary} \label{corollaryentropy}
 Let $\mu$ be a finitely supported aperiodic measure on  $G$  and $f(n)$ be a function such that 
$f(n)  \sqrt{H(n+1)-H(n)} \to 0$
 as $n\to \infty$. 
The distributions of $(G,\mu)$ are individually almost invariant on the scale $f(n)$.

\end{corollary}

{\bf Proof.} 
For  some $C>0$ and any $s$ in the support of $\mu$
$$
\sum_{h}|(\mu^{*n}(h)   -\mu^{*n}(hs)) \le C \sqrt{H(n+1)-H(n)},
$$
(see 2) of Lemma 5.1 in \cite{erschlerkarlsson}. In that lemma it is assumed that $\mu(e)>0$ and this assumption can be easily replaced by aperiodicity).
The Corollary follows therefore from Corollary \ref{corollary1}.


\begin{lemma} \label{lemmaconditionned} [Conditioned measures,  individual almost invariance]
1) Let $\mu_n$ be a sequence of probability measures on $G$, $f_n: G \to \mathbb{R}_+$.
Let $\alpha_n$ be a sequence of events  of for $\mu_n$ such that $\mu_n(\alpha_n)\to 1$ as $n \to \infty$. Let $\mu'_n$ be the conditional measure $(\mu_n| \alpha_n)$. If $f_n$ are $(K,\epsilon)$ are individually almost invariant with respect to $\mu'_n$ and $\Omega_n$, then for any $\epsilon'$ and any sufficiently large $n$ $f_n$ are $ K (1+\epsilon'), \epsilon(1+\epsilon')$ individually almost invariant with respect to  $\mu_n$ and $\Omega_n$.

2) Let $\mu_n=\mu^{*n}$ and let 
 $\alpha_n$ be a sequence of events on $G^n$
  viewed as the space of trajectories of length $n$ for the random walk $(G,\mu)$ such that $P[\alpha_n]\to 1$ as $n \to \infty$. Let $\mu'_n$ be the  projection on the $n$-th coordinate of $G^n$ of the conditional measure $(\mu_n| \alpha_n)$. If $f_n$ are $(K,\epsilon)$  individually almost invariant with respect to $\mu'_n$ and $\Omega_n$, then for any $\epsilon'$ and any sufficiently large $n$ $f_n$ are $ K (1+\epsilon'), \epsilon(1+\epsilon')$ individually almost invariant with respect to  $\mu_n$ and $\Omega_n$.

\end{lemma}

{\bf Proof.}
1) Since $\mu_n(\alpha_n) \to 1$, it holds $ \sum_{g\in G}|\mu_n(g)- \mu'_n(g)| \to 0$. Therefore, 
for any $\epsilon'>0$ and any sufficiently large $n$ $\mu_n(g\in G: |\mu_n(g)/ \mu'_n(g)-1| \ge \epsilon')
\le  \epsilon'$.  Fix $h_n \in \Omega_n$.
Observe also that for any $\epsilon'>0$ and any sufficiently large $n$ $\mu_n(g\in G: |\mu_n(gh_n)/\mu'_n(gh_n)-1| \ge \epsilon') \le \epsilon'$. 
Put $U^h_n= \{g\in G:  |\mu_n(gh)/\mu'_n(gh)-1| \le \epsilon_1   {\mbox  {\rm and }}                 |\mu_n(gh)/\mu'_n(gh)-1| \le \epsilon_1\}$. 
Suppose that $V_n$ is such that $\mu_n(V_n) \ge 1-\epsilon$ and 
$$
1/K \le \mu_n'(g_nh_n)/\mu_n(g_n) \le K
$$ 
for $n \ge N$, all $g_n \in V_n$ and all $h_n\in \Omega_n$.
 Consider $\tilde{V}_n = V_n \cap U_n$.
  We have $\mu_n(\tilde{V}_n) \ge 1- \epsilon- 2\epsilon_1$. Observe that for any $g_n \in \tilde{V}_n$ and any $h_n \in \Omega_n$ it holds
$$
    1/K (1-\epsilon')^2    \le \mu_n(g_nh_n)/\mu_n(g_n) \le K (1+ \epsilon')^2
$$  
  
This implies the statement of 1).

2) Analogously to 1), we observe that $\sum_{x\in G^n} |\mu_{n}(x) -\mu_n''(x)| \to 0$, where $\mu''_n$ is the conditional measure for $\mu_n$ for the condition $\alpha_n$.  This implies that  
for any $\epsilon'>0$ and any sufficiently large $n$  it holds $\mu_n(x\in G^n: |\mu_n(x)/ \mu''_n(x)-1| \ge \epsilon_1)
\le  \epsilon_1$. Therefore, the same condition holds for the projection of $\mu_n$ and $\mu''_n$ to the $n$-th coordinate of $G^n$:
for any $\epsilon'>0$ and any sufficiently large $n$ $\mu_n(g\in G: |\mu_n(g)/ \mu'_n(g)-1| \ge \epsilon_1)
\le  \epsilon_1$. Similarly to the proof of 1), this implies the statement of 2).






\begin{lemma}[A general upper bound for $\bar{r}(n)$]   \label{lemmanotalmostinvariantindividually}
Let $g_n \in G$ and  $A_n$ be an event on the space of trajectories of length $n$. Suppose that
$\mu^{*n} (A_n) >c $ for some $c>0$ and  all sufficiently large $n$. Suppose also  that $\mu^{*n} (g_n A_n)$ tends to zero as $n$ tends to zero. Then the distributions of $(G,\mu)$ are not weakly almost invariant with respect to $g_n$. In particular, the radius of individual almost invariance satisfies $\bar{r}(n) \preceq l(g_n)$.

\end{lemma}

{\bf Proof.} Suppose that the distributions of $(G,\mu)$ are  weakly almost invariant with respect to $g_n$.
Then on some subsets of probability close to one the ratio $\mu^{*n} (g_n x)/ \mu^{*n}(x) >C$, for some $C>0$ and all sufficiently large $n$. This would in particular hold for  subsets of $A_n$ of probability at least $c/2$. This would imply that $\mu^{*n} (g_n A_n) \ge c/2 C$ for all sufficiently large $n$.

In certain situation one gets a better upper bound for $r_{\rm a.i.}(n)$ from using the following:

 \begin{lemma}[A general upper bound for $r_{\rm a.i.}(n)$]   \label{lemmanotalmostinvariant}
 Let $A_n$ be a sequence of events on the space of trajectories of length $n$ for the random walk $(G,\mu)$
 such that $\mu^{*n}(A_n) \to 1$.
 Let $\mu'_n$ be the conditioned measure $\mu^{*n}$ with respect to $A_n$.
 Suppose that for any $V_n \subset G$, with $\mu^{*n}(V_n) \to 1$ and any $K>0$ there exist
$g_n\in G$, $l(g_n) \le r(n)$ and $h_n  \in V_n$ such that $\mu'_n(g_nh_n)/\mu'(g_n) \ge K$.
Then the same conclusion holds for $\mu^{*n}$: for any $V_n$ as above and any $K>0$ there exist 
$g_n\in G$, $l(g_n) \le r(n)$  such that $\mu_n(g_nh_n)/\mu(g_n) \ge K$
In particular, $\mu^{*n}$ are not almost invariant on the scale $r(n)$
(and hence the radius of individual almost invariance satisfies $r_{a.i.}(n) \preceq r(n)$).

\end{lemma}

{\bf Proof.}
Observe that since $\mu^{*n}(A_n) \to 1$ we know that 
 for any $\epsilon>0$ there exist $N$ such that such that for all $n\ge N$ and all  $g\in G$
$$
\mu^{*n}(g)   = P[X_n= g, \mbox{$A_n$}] +  P[X_n= g, \mbox{$A_n$ does not hold}] \ge P[X_n= g, \mbox{$A_n$}] 
  \ge    (1-\epsilon) \mu'_n(g)
$$  

Observe also that there exists $W_n$, $\mu^{*n}(W_n) \to 1$  and a sequence $\epsilon_n$, tending to zero,  such that for all $g\in W_n$
$|\mu^{*n}(g)/  \mu'_n(g)-1| \le \epsilon_n$.
This implies the statement of the lemma.

\subsection{Comparison with the radius of almost constancy in the neighborhood of the identity} \label{fixedpoint}

\begin{lemma} \label{lemmaamenable}
Let $G$ be a finitely generated group and $\mu$ be a symmetric non-degenerate aperiodic
 measure.  The following conditions are equivalent
\begin{enumerate}

\item
 There exists  a sequence
$r(n)$, tending to infinity as $n$ tends to infinity, such that the distributions of $(G,\mu)$ are almost constant in the neighborhood of $e$   of radius $r(n)$.

\item
 There exists  a sequence
$r(n)$, tending to infinity as $n$ tends to infinity, such that the distributions of $(G,\mu)$ are weakly almost constant in the neighborhood of $e$ of radius $r(n)$.

\item  $G$ is amenable.

\end{enumerate}
\end{lemma}

{\bf Proof.}
First show that $(2)$ implies $(3)$. Suppose that $G$ is non-amenable.
There exists $c<0$ such that for all $n>N$ it holds $\mu^{n+1}(e) \le c \mu^{*n}(e)$.

This implies that  for any $K>0$ there exists $m$ such that  $\mu^{*(n+m)}(e) \le  \mu^{*n}(e)/K$ for all $m$ and all sufficiently large $n$.
Observe that 
$\mu^{*(n+m)}(e) = \sum_{g\in G} \mu^{*n}(g)\mu^{*m}(g^{-1})$.

Consider $M$ such that $\mu^{*m}(B(e,M)) \ge 1/2$.
Since $\mu$ is symmetric, we know that $\mu^{*n}(g) \le \mu^{*n}(e)$ for all even $n$ and all $g$.

Observe that for all sufficiently large $n$ there exists an element $h_n$ of length at most $M$ such that
$\mu^{*n}(h_n) \le 2/K \mu^{*n}(e)$.

Now we show that $(3)$ implies $(1)$. 
The ratio theorem of Avez  (\cite{avezratio}) states that for any symmetric non-degenerate measure $\mu$ on an amenable group $G$, such that $\mu^{*m}(e) \ne 0$ for some odd $m$ and any $g \in G$ it holds
$$
\lim_{n \to \infty} \mu^{*n}(g)/\mu^{*n}(e) =1.
$$

Chose $N(i)$ such that for any $n>N(i)$ and any $g$ of length at most $i$ 
$$
1-1/i < \mu^{*n}(g)/\mu^{*n}(e) < 1+1/i.
$$
Let $r(n)$ be the maximum of $i$ such that $N(i) <r(n)$.
It is clear that for any $g_n$ such that $l(g_n) \le r(n)$ it holds
$\mu^{*n}(g_n)/\mu^{*n}(e)  \to 1$.

Finally, it is clear that $1$ implies $2$.

\begin{example} [Revelle, Theorem 1 \cite{revelle}] \label{examplerevelle}
Let $G=\mathbb{Z}\wr F$, $F$ is a finite group and $\mu$ be an aperiodic  switch-walk-switch finitely supported measure on $G$. For any $\epsilon>0$ there exists $c>0$ such that for all sufficiently large $n$ and all $g$ of length at most $cn^{1/3}$ it holds $(1-\epsilon)\le \mu^{*n}(g)/\mu^{*n}(e) \le 1$.
For any $\alpha(n): \alpha(n)/n^{1/3} \to \infty$  there exists $g_n \in G$, $l(g_n)\ge \alpha(n)$ such that $\mu^{*n}(g_n)/\mu^{*n}(e)  \to 0$.

\end{example}

Revelle uses the fact that the random walk is switch-walk-switch in the proof of Theorem 1 in \cite{revelle} to obtain more precise description of the decay of transition probabilities in the neighborhood of the identity.
 The statement of Example \ref{examplerevelle} can be proven for any finitely supported aperiodic random walk on $\mathbb{Z}\wr F$

\begin{remark} Let $G$ be an amenable finitely generated group and $\mu$ be a non-degenerate symmetric finitely supported measure such that the Poisson boundary of $(G,\mu)$ is non-trivial. The radius of almost constancy of transition probabilities  of $(G,
mu)$i s not trivial (that is,  $r_{\rm a.i.}(n)$ tends to infinity), while the radius of almost invariance is trivial ($R(n)$ is bounded). 
\end{remark}

While in Abelian groups $R(n)$ and $r_{\rm a.i.}(n)$, as well as $\bar{r}(n)$ are asymptotically equivalent,
the remark shows that $r_{\rm  a.i.}(n)$ could be asymptotically smaller  than $R(n)$ (and not asymptotically equivalent). Theorem \ref{theorem1} shows that for 
$G=\Z\wr F$  we have $r_{\rm a.i.}(n)\sim n^{1/2}$, and since $R(n) \sim n^{1/3}$ we see that 
$r_{\rm a.i.}(n)$ can be larger than $R(n)$.

\section{Proof of Theorem \ref{theorem1}} \label{prooftheorem1}

\subsection{Auxiliary facts about random walks on $\mathbb{\Z}$} \label{auxiliary}

For a random walk on  $\mathbb{Z}$ we denote by $Min_n$ and $Max_n$ the minimal and the maximal point of $\mathbb{Z}$, visited at least once until the moment $n$.

\begin{lemma} \label{firstlemmaline}

Consider a random walk on $\mathbb{Z}$ defined by a symmetric finitely supported measure $\mu$ such that the support of $\mu$ is not equal to the identity. 
For any $\epsilon>0$ there exists $c>0$ such that with probability $\ge 1-\epsilon$ with respect to $\mu^{*n}$ it holds
$$
 -c \sqrt{n} \le {\rm Min_n} ,     {\rm Max_n} \ge c \sqrt{n}
$$

\end{lemma}

{\bf Proof.} 
Recall that by a result of Erd{\"o}s and Kac \cite{erdosKac} the limit as $n\to \infty$ of the probability
$P[{\rm Max_n} \le c \sqrt{n}] = C \int_{x=0}^c \exp(-t^2/2) \partial t$,
for any non-degenerate $\mu$ on $\mathbb{Z}$ with zero mean and finite second moment.
Therefore,  for  $\epsilon>0$ there exists $c>0$ such that with probability $\ge 1-\epsilon$ with respect to $\mu^{*n}$ it holds
$ {\rm Max_n} \ge c \sqrt{n}$. Since ${\rm Min_n}$ is equal to ${\rm Max_n}$ for the trajectory, reflected at $0$, it follows that for any $\epsilon>0$ there exists $c>0$ such that with probability $\ge 1-\epsilon$ with respect to $\mu^{*n}$  it holds $-c \sqrt{n} \le {\rm Min_n}$.

For a random walk on $\mathbb{Z}$ denote by $I_n$ the interval of $\mathbb{Z}$ visited until time instant $n$. Observe that if the support of $\mu$ belongs to $\{ -1,0,1 \}$, then this interval coincides with the sites visited until the time instant $n$.

\begin{lemma} \label{lemmaray1}
Consider a simple random walk on $\Z$, reflected at $0$. 
For any $C, \epsilon>0$ ther exist $a>0$ such that $n$ step transition probabilities satisfy
$$
| p_n(0, x+i)/ p_n(0, x) -1    | \le \epsilon
$$

for any $x: 0\le x \le C \sqrt{n}$ and any $i: 0\le i \le a \sqrt{n}$.

\end{lemma}

{\bf Proof.}

Observe that $\exp(- K x^2/n) / \exp(-K (x+i)^2/n)$ is close to one whenever $ 0\le x \le C \sqrt{n}$ and any $0\le i \le a \sqrt{n}$ and  taht
$p_n(0,0)/p_n(0,x)$ is close to $\exp(- K x^2/n)$ whenever   $ 0\le x \le C \sqrt{n}$ (as it follows from the Local Limit theorem, see e.g.   Theorem  2.3.11 in \cite{lawlerlimic}).






\subsection{Proof of Theorem \ref{theorem1}}

Take an element $z \in G$, $h=(a,f)$, $a\in \mathbb{Z}$, $f$ is a finitely supported function on $\Z$ with values in $F$.
For such function, we denote by $I(h)=I(f)$ the minimal interval, containing the support of $f$.
Denote by $J(h)$, the minimal interval containing $a$, $e$ and $I(f)$.

\begin{lemma} \label{lemmasupportcontainsinterval} Let $\mu$ be a symmetric non-degenerated finitely supported measure on $G= \mathbb{Z} \wr F$, $F$ is a finite group containing at least $2$ elements. Denote by $X_1$, $X_2$, \dots, $X_n$ a trajectory of a random walk $(G,\mu)$, $X_i= (x_i, f_i)$, $x_i \in \mathbb{A}$, $f_i: \mathbb{Z} \to F$.

1) For any $\epsilon>0$ there exists $ C>0$ such that 
$I(X_n)$  belong to the interval $[-C\sqrt{n} , C\sqrt{n}]$ with probability greater than $1-\epsilon$.

2) For any $\epsilon>0$ there exist $ \epsilon_1>0$ such that $I(f)$ contains  the interval $[-\epsilon_1 \sqrt{n} , \epsilon_1 \sqrt{n}]$ with probability greater than $1-\epsilon$.

3) For any $\epsilon>0$ there exist $ \epsilon_1>0$ such that $I(f)$ contains  the interval $[x_i-\epsilon_1 \sqrt{n} , x_i+ \epsilon_1 \sqrt{n}]$   with probability greater than $1-\epsilon$.

\end{lemma}

{\bf Proof.}
1) Follows from the fact that the support of $g$ belongs to the set, visited by the projection of the random walk to $\Z$ up to the moment $n$ and from the fact, that any finitely supported symmetric random walk on $\Z$, stays with a probability close to one (with respect to $n$ step distribution) in a ball of radius $C \sqrt{n}$.

2)  Lemma \ref{firstlemmaline}  implies 
 that with the probability close to one, with respect to $n$ step distribution, the projection of the random walk on $\Z$ has visited some points $X$ and $Y$,  $X \le -\epsilon_2 \sqrt{n}$ and $ Y\ge \epsilon_2 \sqrt{n}$.
Let $S$ be a finite generating set of $G$ and $C_S$ be the maximum of $l_S(s)$, where $s \in \sup \mu$.

Suppose that the random walk has visited $X,Y$ as above.
Observe that the conditional probability that  $I(X_n)$ contains some elements $-L$ and $M$, with $L,M \ge 1/2 \epsilon_2 \sqrt{n}$ is greater than
$\exp(-C/C_S \epsilon_2  \sqrt{n})$, for some $C>0$, all $\epsilon_2$ and all $n$.

This probability is close to one if $\epsilon_2$ is small enough, and 
we have proved  2) the lemma.

Now let $K_n$ be the maximal point of $\mathbb{Z}$ visited by the projection of the random walk to $\mathbb{Z}$ until the instant $n$, and let
$K'_n$ be the minimal point of $\mathbb{Z}$ visited by the projection of the random walk to $\mathbb{Z}$ until the instant $n$.
 Observe that $K_n-x_n$ (as well as $x_n-K'_n$)
  has the same distribution as a random walk on $\mathbb{Z}$, reflected at $0$. This implies that $K_n-x_n \ge 2 \epsilon_1 \sqrt{n}$ and that 
 $x_n-K'_N \ge 2 \epsilon_1 \sqrt{n}$  with probability close to one, for all sufficiently large $n$.
 Observe  that at least one of the integer points in the interval $[K_n -\epsilon_1 \sqrt{n}, K_n]$ belongs to the support of $f_i$, with probability close to one.
 With the same argument we conclude that at least one of the integer points in the interval $[K_n', K'_n +\epsilon_1 \sqrt{n}]$ belongs to the support of $f_i$, with probability close to one.
This implies 3) and completes the proof of the Lemma. 

Now we start by proving the second claim of the theorem.

\begin{lemma} \label{lemmaswitch}
Let $\mu$ be a measure defining a standard simple switch-walk-switch random walk on $G=\mathbb{Z}\wr F$. 

1)Suppose that  $g_0=(e,h)$  and  that $z\in G$ is such that 
 $J(g_0) \subset J(z)$.
Then for any $n>0$ $\mu^{*n}(z)= \mu^{*n} (g_0z) $. 

2) Suppose that  $g_0=(e,h)$  and  that $z =(x,h')\in G$ is such that 
$x J(g_0) \subset J(z)$.
Then for any $n>0$ $\mu^{*n}(z)= \mu^{*n} (zg_0) $. 

\end{lemma}

{\bf Proof}.
Take some  $z =(a,f) \in \mathbb{Z}\wr F$. By assumption on the measure, the support of its projection to $\Z$ is equal to $\{-1,0,1\}$.
Observe that if the projection of the random walk $(G,\mu)$ visits $X,Y \in \mathbb{Z}$ until some time instant $n$, then this projection visits all integer points of the interval $[X,Y]$.
Let $X_n=(a_n, f_n)$ be the trajectory of the random walk on $G$.
Observe that if we condition the random walk to the condition that the projection to $\mathbb{Z}$ visits some interval $I$, then for each integer $k\in I$ the value of $f_n(k)$ takes all values in $F$, and that these values are independent for all $k$ (in this interval)

This implies that the transition probability
$$
\mu^{*n}(z) = \sum_{I} P[I(X_n)=I, I(z) \subset I, \mbox{ and } P(X_n)=a ] \frac{1}{\# F^\# I} =
$$

$$
=\sum_{I: I(z) \subset I} P[I(X_n)=I \mbox{ and } P(X_n)=a ] \frac{1}{\# F^\# I}
$$

This implies the statement 1) and 2) of Lemma \ref{lemmaswitch} and concludes the proof of the second claim of the theorem.

Now we prove the first claim of the theorem.

Let us show that for each $\epsilon >0$ there exists $a, C>0$ and $V_{n,\epsilon} \subset G$, with 
 $\mu^{*n}(V_{n,\epsilon}) \ge 1-\epsilon$ for all sufficiently large $n$  such that the following hold:
 
 1)  $\mu^{*n}(g)  \ge \exp(- C \sqrt{n})$ for all sufficiently large $n$ and any $g \in V_{n,\epsilon}$.
 
2)  For all sufficiently large $n$,  and all $g =(x,f) \in V_{n,\epsilon}$ the interval $I(f)$ contains $[x-2a \sqrt{n}, x + 2 a \sqrt{n}]$.


Indeed, observe the probability (with respect to $\mu^{*n}$) of the set of $z$ satisfying 1) is close to one, as follows from Lemma \ref{lemmaswitch} combined with 
 the fact that the projection of the random walk to $\mathbb{Z}$ 
stays inside the interval $[-B' \sqrt{n}, B' \sqrt{n}]$ with positive probability, which is close to one if $B'$ is large enough

Observe also the probability (with respect to $\mu^{*n}$) of the set of $z$ satisfying 2) is close to one in view of statement 3) of Lemma \ref{lemmasupportcontainsinterval}.

Now we estimate the ratio $\mu^{*n}(gx)/\mu^{*n}(g)$ for $g \in V_{n,\epsilon}$.

Let $z$ be the generator of $\mathbb{Z} \subset \mathbb{Z} \wr F$.
Observe that the transition probability $\mu^{*n}(g')$, conditioned to the fact that the projection of the random walk on $\mathbb{Z}$ has visited an interval of length $L$ is less or equal to $\exp(-C'L)$ for some $C'>0$ and all $l$. In view of Lemma \ref{lemmaray1} this shows that 
$$
|\mu^{*n}(g z^i)/\mu^{*n}(g)-1| \le \epsilon ,
$$
for all $g \in  V_{n,\epsilon}$ and all $i \le a \sqrt{n}$, whenever $a$ is a small enough constant.

Now observe that 2) implies that 
 for all sufficiently large $n$, all  $i: -a  \sqrt{n} \le i \le a \sqrt{n}$ and all $g \in V_{n,\epsilon}$ 
$ J(g z^i)$ contains a ball of radius $a \sqrt{n}$, centered at the projection of $g z^i$ to $\Z$.

Combining this observation with the second statement of Lemma  
we see  that $\mu^{*n}(g z^if) =\mu^{*n}(g z^i)$ for all $f$ which has the trivial projection to $\Z$ and such that $J(f) \subset [-a\sqrt{n}, \sqrt{n}]$.

Observe that any element $g$ of the wreath product can be written as a product $g=x_0 g_0 $, where $g_0$ has trivial projection to $\Z$ and $x_0 \in \Z$. It is clear that for any word metric $l_S$ on the wreath product there exists a constant $C>0$ such that for all $g$ the length of the elements $g_0$ and $x_0$ in their decomposition above satisfy $l_s(g_0), l_s(x_0) \le C l(g)$.

We see therefore that  
$$
|\mu^{*n}(g h)/\mu^{*n}(g)-1| \le \epsilon 
$$
for any $h$ of  word length at most $Ca \sqrt{n}$.

This completes the proof of the theorem.

\section{More on wreath products and iterated wreath products} \label{morewreath}

Let $C$ be a wreath product of $A$ and $B$, and $c
\in C$, $c=(a, f)$, where $f: A \to B$.  Given $a \in A$, denote by $c((a))$ the value $f(a)$. Take  $C= A\wr( A \wr B)$ and $a_1, a_2 \in A$. Observe that $c((a_1))$ is an element of $A \wr B$, denote by 
   $c((a_1,a_2))$ its value at $a_2$. Likewise, for $i$ times iterated wreath product $C \wr (A \wr(A \wr  \dots ) \wr B ) \dots )$ we consider $c((a_1, \dots , a_i))$ defined inductively as the value of $c((a_1, \dots , a_{i-1}))$
at $a_i$. If $a_1 =a_2 \dots =a_i e_A$, we also use the notation
$c(e,i)$ for $c((a_1, \dots, a_i))= c((e_A, \dots, e_A))$.

\begin{lemma} \label{lemmaiterated} Let $G_i (H)$ be $i$ times iterated wreath product of $\Z$ with $H$ ($\Z$ is the group that acts).
Consider switch-walk-switch symmetric finitely supported random walk  $\mu_i$ on $G_i$ and suppose that the drift function of the corresponding random walk on $H$ is $\sqrt{n}$.
Then the expected value after $n$ steps of the random walk $(G_i, \mu_i)$ 
$$
E[l(X_i((e,i)] = n^{1/2^{i+1}}.
$$

\end{lemma}

{\bf Proof.}  Let $L_0(n)$ be the number of  visits of $e$ after $n$ steps of  a random walk on $\Z$.
The expectation of $L_0(n)$ is $\sim C \sqrt{n}$.
 Consider the conditional event that $L_0=\bar{L}_0$. If we condition moreover by the  fact that the point of $\Z$ visited by the random walk  at the moment $n$ is not equal to $0$, then 
  the value of $X_i(e)$ is obtained as a position after $2\bar{L}_0(n)-1$  steps of a switch-walk-switch symmetric finitely supported switch-walk-switch random walk on  $G_{i-1}(H)$.
Otherwise, if we condition  by the  fact that the point of $\Z$ visited by the random walk  at the moment $n$ is  equal to $0$, then 
  the value of $X_i(e)$ is obtained as a position after $2\bar{L}_0(n)$  steps of the above mentioned switch-walk-switch symmetric finitely supported switch-walk-switch random walk on  $G_{i-1}(H)$.

Observe that $\sqrt{n}$ is concave, and therefore the expectation of $\sqrt{L_0(n)}$ is not greater than $\sqrt{E[L_0(n)]}$, and this implies that $E[l(X_1(e,2)] \le C n^{1/4}$. Observe that there exists $a, p>0$ such that with positive probability ($\ge p$) it holds $L_0(n) \ge a \sqrt{n}$.
 This implies that $E[l(X_1(e,2)] \ge p K n^{1/4}$. Therefore,   
$E[l(X_1(e,2)] \sim n^{1/4}$ and we have proved the statement of the lemma for $i=1$. 

Arguing by induction on $i$ and using concavity of $n^{1/2^i}$, we obtain the claim of the lemma for any $i\ge 1$.

\begin{corollary} \label{coriterated}
 Let $G_i (H)$ be $i$ times iterated wreath product of $\Z$ with $H$ ($\Z$ that acts each time).
Consider switch-walk-switch symmetric finitely supported random walk  $\mu_i$ on $G_i$ and suppose that the drift function of the corresponding random walk on $H$ is $\sqrt{n}$.
Then the radius of individual almost invariance  $\bar{r}(G_i,\mu_i)(n) \le C / n^{1/2^{i+1}}$.

\end{corollary}

For example, Corollary \ref{coriterated} can be applied to $H =\Z$ and states that simple switch-walk-switch random walks on $\Z\wr Z$ have the radius of individual almost invariance  asymptotically not larger than $n^{1/4}$, 
that in case of random walks on $\Z \wr \Z \wr \Z $ the corresponding radius of individual almost invariance is at most 
$n^{1/8}$ etc. 

{\bf Proof.} Follows from Lemma \ref{lemmaiterated} and Lemma \ref{lemmanotalmostinvariantindividually}.

Since radius of almost invariance is less or equal to the radius of individual almost invariance,  Corollary \ref{coriterated} implies in particular that
 $r_{\rm a.i.}(G_i,\mu_i)(n) \preceq C / n^{1/2^{i+1}}$.

This upper bound for $r_{\rm a.i.}(G_i,\mu_i)(n)$ can be improved:

\begin{proposition} \label{propositioniteratedwreath}
Let $G_j $ be $j$ times iterated wreath product of $\Z$ with $\Z$.
Consider switch-walk-switch symmetric finitely supported random walk  $\mu_j$ on $G_j$.
Then radius of almost invariance $r_{\rm a.i.}(G_j,\mu_j)(n) \le C  \frac{  n^{1/2^{j+1}}} { \sqrt{\ln(n))}}$. 

\end{proposition}


\begin{corollary} \label{corollary1comparewithdrift}
Under the assumption of  Proposition \ref{propositioniteratedwreath} ,   the radius of almost invariance of this walk is asymptotically strictly smaller than $n/L(n)$.
\end{corollary}

We recall that $L(n)$ denotes the drift function of the random walk $(G,\mu)$.

{\bf Proof.} 
Follows from the fact that
the drift of the standard switch walk switch random walk on $G_i$ is asymptotically equivalent to $n^{1-1/2^{i+1}}$ (\cite{erschlerdrift}. The random walk considered in \cite{erschlerdrift} are "switch or walk" and the same argument works for "switch walk switch" random walks).

The assumption on the simple random walk in Corollary \ref{corollary1comparewithdrift} is not important.

{\bf Proof of Proposition \ref{propositioniteratedwreath}.} Consider $G_1 = \Z \wr \Z$ and a symmetric finitely supported switch-walk-switch random walk  $\mu_1$ on $G_1$.
Let $X_i= (z_i, f_i)$ be the trajectory of the random walk $(G_1, \mu_1)$, where $z_i \in \Z$ and $f_i$ are finitely supported functions from $\Z$ to $\Z$.

Take $a: 0< a< 1/2$. Take some $n\ge 1$. We want to show that  that, for some $C>0$,  with probability close to one  with respect to $\mu_1^{*n}$ there exists $y \in \Z$ such that $f_n(y) \ge C \sqrt{n} \sqrt{\ln(n)}$.
Indeed, let $t_{i,n}'$ be the number of visits of $i$ by the projection $z_i$ of the random walk to $\Z$ until the moment $n$, and let $t_{i,n} = 2t_{i,n}$ if $i \ne 0, i \ne z_i$; $t_{0,n} = 2t'_{0,n}+1$,
$t_{i,n}= 2t'_{0,n}+1$ for $i=z_n$ if $z_n \ne 0$ and $t_{0,n} = 2t'_{0,n}$ if $z_n=0$. Observe that the value $f_n(i)$ is obtained as a position of after $t_{i,n}$ steps of the random walk on $\Z$, and these random walks are independent for $i \in \Z$.
Normalized $t_{i,n}'$ converge to the local time of Brownian motion,
and hence with positive probability (close to one if $\epsilon$ is small enough) it is true that $t_{i,n}' \ge  \epsilon \sqrt{n}$ for all $i \le \epsilon n$, and in particular for all $i \le n^a$.
Therefore,  $t_{i,n} \ge  \epsilon \sqrt{n}$  for all $i \le n^a$.

\begin{lemma} \label{lemmaestimatebig}
Let $Y_i$ be a trajectory of a symmetric finitely supported random walk on $\Z$ (with the support of the defining measure that contains at least two points).
There exist $C_1, C_2>0$ such that for any $n$ and any $L\ge \sqrt{n}$
the probability that $Y_n \ge L$ is greater or equal to $C_1 (\sqrt{n}/L) \exp(- C_2  (L^2/n))$
\end{lemma}

{\bf Proof.}
Denote by $\mu$ the measure that defines our random walk on $\Z$.
Using Gaussian estimates for the random walk $(\Z,\mu)$ we observe that
probability that $Y_n \ge L$ satisfies
$$
P[Y_n\ge L]= \sum_{i=L}^{\infty} \mu^{*n}(i) \ge \sum_{i=L}^{[L+n/L]} \mu^{*n}(i) \ge \sum_{i=L}^{[L+n/L]} 1/\sqrt{i}\exp(-Ki^2/n)
$$
Since $L\ge \sqrt{n}$, we
observe that for any $i: L \le i \le L+n/L$ it holds $\sqrt{i} \le \sqrt{2L}$ and $i^2/n \le (L+n/L)^2/n \le (2L^2)/n+2$. This implies that that for some $C_1, C_2>0$
$$
P[Y_n\ge L]  \ge C_1 (n/L) (1/\sqrt{n}) \exp(-C_2 L^2/n)  = C_1 (\sqrt{n}/L) \exp(-C_2 L^2/n), 
$$
and this completes the proof of the Lemma.

Lemma \ref{lemmaestimatebig} implies that 
for a symmetric finitely supported random walk  $Y_i$ on $\Z$  it holds $Y_i \ge \sqrt{i} \sqrt{\ln (i)}$ with probability $\ge  1/\sqrt{\ln(i)} \exp(-C \ln i)= (\sqrt{\ln(i)} i^{-C})$. 

Fix some positive $a$ such that $a<1/2$.
Chosing positive $C$ small enough we assure that  $(\sqrt{\ln(i)} i^{-C} < i^a$. In this case  we can consider $n^a$ independent random walks,  making $t_j\ge c n^{1/2}$ steps, at least one of them stays at the time $t_j$ at some point 
$\ge C \sqrt{n} \sqrt{\ln(n)}$.  Therefore,  with positive probability (which is close to one so far as $C$ is small enough) there exist $y : 0< y < n^a< n^{1/2}$ such that $f_n(y) \ge C \sqrt{n} \sqrt{\ln(n)}$.

\begin{lemma} \label{lemmadecayinz}
Consider a symmetric finitely supported random walk on $\Z$ (with the support of the defining measure that contains at least two points).  
For any $\epsilon>0$ there exists $K>0$ such that for  any sufficiently large $n$ and 
any $:T: \sqrt{n} \le T \le n^2/3$
 it holds $\mu^{*n}(T + [K n/T]) \le \epsilon \mu^{*n}(T)$.
\end{lemma}

{\bf Proof.} 

There exists $C_0, C_1, C_2>0$ such that for any $k\le n^{2/3}$
$$
  C_0 1/\sqrt{n} \exp(-C_2 k^2/n) )      \le \mu^{2n}(2k)  \le C_1 1/\sqrt{n} \exp(-C_2 k^2/n) )
$$
(see e.g.  Theorem  2.3.11 and the remark after this theorem in \cite{lawlerlimic}). 
Therefore, for any $k_1, k_2 \le n^{2}/3$ the ratio $\mu{2n}(2k_1)/\mu^{*2n}(2k_2)$ is close up to multiplicative constant to $\exp(-C_2 (k_1^2-k_2^2)/n)$.

This implies the statement of the Lemma.

Now we return to the proof of the proposition.
Fix some positive $a$ such that $a<1/2$.
Take any function $C(n)$ tending to $\infty$.
Consider $g_n=(z,f) \in \Z \wr \Z$ and assume that for some $y$ such that $|y-z|<n^a$,  it holds $f(y)  \ge C_1 \sqrt{n} \sqrt{\ln(n)}$. Observe that for all  $K>0$ there exist $C>0$, depending
 on $K$ and $C$,  such  that the transition probabilities of the projection random walk $\nu$ on $\Z$
satisfy $\nu^{*n}(x- C \sqrt{n}/\sqrt{\ln(n)})/\nu^{*n}(x) \ge K$ for any  positive integer $x: l(x) \ge C_1 \sqrt{n} \sqrt{\ln(n)}$.
 Consider $f': \Z \to \Z$ such that $f'(y)=-C(n) \sqrt{n}/\sqrt{\ln(n)}$ and $f'(x)=0$ for all $x \ne y$.
Consider $h_n\in \Z\wr \Z$ such that $h_n= (0,f')$. Observe that the length of $h_n$ with respect to the word metric on the standard generators of $\Z \wr \Z$ is at most $n^{a}+ C (n)\sqrt{n}/\sqrt{\ln(n)} \le 2 C(n)\sqrt{n}/\sqrt{\ln(n)}$.
Take a   sequence $\delta_n$, tending to $0$,  such that $\delta(n)/C(n) \to \infty$.
  and let $A_n$ be the event on the space of trajectories of length $n$ of $(G,\mu)$ for which the number of visits of $y$ by the projected random walk on $\Z$ is greater or equal than
$\delta_i \sqrt{n}$. Let $\mu'_n$ be the conditional measure of $\mu^{*n}$ with respect to $A_n$.

Observe also that Lemma \ref{lemmadecayinz} implies that $\mu'_n(g_nh_n)/\mu'_n(g_n) \ge  K$. 
Taking in account Lemma \ref{lemmanotalmostinvariant}, this implies that the random walk $(G,\mu_1)$ is not almost invariant on the scale  $C(n) n^{1/4}/ \sqrt{\ln(n)} $ (for none of the the function $C(n)$ tending to infinity). 

Therefore,  the almost invariance radius for the random walk on $G_1=\Z \wr \Z$ satisfies
 $r(G_1,\mu_1)(n) \le C / ( n^{1/4}/\sqrt{\ln(n)}  )$, and this completes the proof of Proposition for  $G_1=\Z \wr \Z$.

To prove the proposition for $G_j$ for any $j \ge 1$, consider the value  $\phi_{j,n}=f_n(\underbrace{e, e, \dots, e}_{j-1}) \in \Z \wr \Z$.  
Observe that $f_n$ is obtained as a value    of a finitely supported symmetric a random walk  $W_n=(z_n, f_n)$ on $\Z \wr \Z$, obtained after $\tau_n$ steps, where $\tau_n$ is the number of visits of $(\underbrace{e, e, \dots, e}_{j-1})$.
Arguing by induction on $j$, we see that the expectation of $\tau_n$ is $C n^{1/2^{j-1}}$ for some $C>0$. Fix some positive $a$ such that $a < n^{1/2^{j+1}}$.
 Analogously to the case $j=1$ we conclude that with probability close to one, for any sufficiently large $n$, there exists
  some integer $y$ such that  $|y-x_n|<  n^a$ and such that $\phi_{j,n}(y)= f_n( \underbrace{e, e, \dots, e}_{j-1},  y) \ge C'  n^{1/2^{j+1}} \sqrt{\ln(n)}$.

Take $C(n)$ tending to infinity and some $c'>0$. Consider $V_n \subset G$ such that $\mu^{*n}(V_n) >c'$. Observe that there exist $g_n =(x_n,f_n) \in V_n$ and
$y_n$ such that  $|y-x_n|<  n^a$ and such that $f_n( \underbrace{e, e, \dots, e}_{j-1},  y) \ge C'  n^{1/2^{j+1}} \sqrt{\ln(n)}$

Now consider a $h_n = h_{n,j} \in G_j$, defined recursively by $h_{n,1}= (e, f')$,  where $f'_1(y_n)=f'(y_n)=[-C(n) n^{1/2^{j+1}}/\sqrt{\ln(n)}]$ and $f'(x)=0$ for all $x \ne y_n$; and $h_{n,k+1}= (e, f'_{k+1})$, where 
$f'_{k+1}(0)=f'_k$ and $f'(x)=0$ for all $x \ne 0$.
It is clear that the length of $h_{n,j}$ is at most $2 C(n)n^{1/2^{j+1}} /\sqrt{\ln(n)}$.
Observe that for any $K>0$ and any sufficiently large $n$ it holds
 $\mu_j^{*n}(g_n h_n)/\mu^{*n}(g_n) \ge K$.
 


This shows that the random walk $(G_j,\mu_j)$ is not almost invariant on the scale  $C(n) n^{1/2^{j+1}}/ \sqrt{\ln(n)}$ (for none of the the function $C(n)$ tending to infinity).
This implies that  the almost invariance radius for the random walk $(G_j, \mu_j)$ is not greater than      $n^{1/2^{j+1}} /\sqrt{\ln(n)}$, and completes the proof of the proposition.

\subsection{Wreath producst of $\Z^2$ with finite groups}

Given a random walk on $(\Z^2,\nu)$, denote by $r_{\rm cov}(n)= r_{\rm cov}(\nu, n)$ the function such that for any $\epsilon$ there exists $c$ with the following property. For all sufficiently large $n$ the set of points of $\Z^2$ visited at least once until the moment $n$  contains the ball of radius $cr_{\rm cov}(n)$ with probability at least $1-\epsilon$.

\begin{proposition} \label{propositionwrz2}
Let $F$ be a finite group containing at least two elements.
Let $\mu$ be a measure on $G=\Z^2 \wr F$ which defines a standard aperiodic switch-walk-switch random walk. Denote by $\nu$ the projection of $\mu$ to $\Z^2$. We assume that the support of $\nu$ consists of $5$ elements: 2 standard generators of $\Z^2$, their inverses and the neutral element.
Then the radius of individual almost invariance of $(G,\mu)$ satisfies $\bar{r}(n)  \succeq r_{\rm cov}(\nu, n)$.

\end{proposition}

{\bf Proof.}
Let $S$ denote some finite generating set of $G$.

We  want to show that for any $r(n): r(n)/r_{\rm cov}(\nu, n) \to 0$ there exists $V_n \subset G$ such that $\mu^{*n}(V_n) \to 1$ and such that for any $ g_n \in V_n$ and 
any $h_n \in G$ such that $l_S(h_n) \le r(n)$ it holds $\mu^{*n}(g_nh_n)/\mu^{*n}(g_n) \to 1$ as $n\to \infty$. 



Analogously to the proof of Theorem \ref{theorem1}, it is sufficient to prove this statement separately for $h_n$ having  trivial projection to $\Z^2$ ($h_n$ belongs to $\sum_{\Z^2} F$) and $h_n \in Z^2$.

First observe that $r_{\rm cov}(\nu, n) \le \sqrt{n}$ since the number of points, visited until the moment $n$ is at most $n$.

Take any function $r'(n)$ such that $r(n)/ r'_(n)$ tends to $0$. Put $r''(n)= \sqrt{r_{\rm cov}(n)r'(n)}$.
We have $r(n) \le r'(n) \le r''(n) \le r_{\rm cov}$ and $(r')^2/n, (r'')^2/n \to 0$ as $n \to \infty$.
Observe that $r''(n)/r_{\rm cov}(n) \to 0$ and $r'(n)/r''(n) \to 0$.
Denote by $\pi$ the quotient map from $G$ to $\Z^2$.

Consider the event $\alpha_n$ on the space of trajectories of length $n$ for the random walk $(G,\mu)$ that consists of the trajectories $X_1, X_2, \dots X_n$, $X_i=(x_i,f_i)$ such  that

(1) the projection of the trajectory to $\Z^2$ visits all points of the ball of radius  $2 r''(n)$, centered at the origin, until the moment $n$;

(2) the projection of the trajectory to $\Z^2$ visits all points of the ball of radius  $2 r''(n)$, centered at $x_n$, until the moment $n$;

(3) Moreover,  put $T_n=   (r'(n))^2$ . The projection of the trajectory to $\Z^2$ visits all points of the ball of radius  $2r''(n-T_n)$, centered at $x_{n-T_n}$, until the moment $n- T_n$ and
$|x_{n-T_n}, x_n| \le r''(n)$.

Observe that (3) implies in particular, that the trajectory of the random walk visits until the moment $n-T_n$ all the points in the ball of radius $r''(n)$, centered at $X_n$.

Denote by $\mu_n$ the conditional measure $\mu^{*n} | \alpha_n$.

First let us show probability of $\alpha_n$ tends to $1$ as $n \to \infty$.
Indeed, the property (1) holds with probability close to one, as it follows from the definition of $r_{\rm cov}$ and from the fact that  $2 r''(n)/r_{\rm cov}(n) \to 0$  as $n \to \infty$.
To prove that (2) holds with probability close to one, consider the "inverted trajectory" of length $n$: $e$, $(X_n X_{(n-1)}^{-1} )^{-1}$,   $(X_n X_{(n-1)}^{-1} )^{-1} (X_{n-1}X_{(n-2)}^{-1} )^{-1}$ \dots and observe that the trajectory $e, X_1, \dots X(n)$ visits the ball of radius $f(n)$ centered at $X_n$ if and only the inverted trajectory visits the ball of the same radius, centered at the origin.
The first property of (3) holds with probability close to one, as follows from (2) applied to $n' = n-T_n$. The second property of (3) holds since $r''(n)/\sqrt{T_n} \to \infty$ as $n \to \infty$.


Now take $h_n$ such that $l_S(h_n) \le r(n)$. First suppose that $h_n \subset \sum_{\Z^2} F$.
Using $r'(n) \le r''(n)$ and  (1) in the definition of $\alpha_n$  we
observe that for any $g \in G$ $\mu_n(h_n g |\alpha_n) = \mu_n (g|\alpha_n)$.
Since $\mu_n(\alpha_n) \to 1$ as $n \to \infty$, this implies that for some sets $V_n$ such that $\mu_n(V_n)\to 1$ it holds
$\mu_n(h_n g )/\mu_n (g)$ is close to one for any $g \in V_n$ and any $h_n$ such that $l_S(h_n) \le r'(n)$. Since $\mu$ is symmetric, this also implies that for some $V_n$ such that  $\mu_n(V_n)n\to 1$,
any $g \in V_n$ and any $h_n$ such that $l_S(h_n) \le r'(n)$ it holds $\mu_n( g h_n )/\mu_n (g)$ is close to one.

Now suppose that $h_n \in \Z^2$, $l(h_n) \le r(n)$. For the projection of the random walk to $\Z^2$, consider the conditional event corresponding to the properties 1) -3).  Let $\nu_n$ be the corresponding conditional measure. Observe that that for $n$ tending to infinity
$$
P_{\nu_n}[x_{n-T_n} =x',  x_n = x'']/ P_{\nu_n}[x_{n-T_n} =x',  x_n = x''h_n] \to 1
$$

Therefore, for some set $V_n$ such that $\mu^{*n}(V_n) \to 1$, any sufficiently large $n$ and any $g\in V_n$ the ratio $\mu_n( g h_n )/\mu_n (g)$  is close to one.

Now we apply 2) of Lemma \ref{lemmaconditionned}.
We observe that for  any $h_n \in G$ such that $l(h_n) \le r(n)$ there exist $V_n$ such that $\mu^{*n}(V_n) \to 1$, such that for 
 any $g_n \in V_n$ it holds $\mu^{*n}( g _n h_n )/\mu^{*n} (g) \to 1$ as $n$ tends to $\infty$.

This completes the proof of Proposition \ref{propositionwrz2}.

\begin{corollary} \label{corollary2lowerboundz2}
Under the assumption of Proposition \ref{propositionwrz2}, the radius of individual almost invariance  of the simple random walk $(G,\mu)$ is asymptotically strictly  larger than $n/L(n)$. Here $L(n)$ denotes as before the drift function of the random walk $(G,\mu)$.

\end{corollary}

{\bf Proof.} It is sufficient to use the estimate of covered balls,  centered at the origin in $\Z^2$  due to R{\'e}v{\'e}sz (see e.g . Theorem 24.2 \cite{revesz}, Chapter 24) that states that for some $C>0$ 
$$
\liminf_{n \to \infty} P[ (\ln (R_{\rm cov}(n))^2/\log(n) >z] \ge \exp(- C z)
$$

(A precise limit distribution (of the constant above) for  size of the ball in $\Z^2$, centered at the origin, and covered by the trajectory of length $n$ is due to Dembo, Peres, Rosen and Zeitouni (\cite{demboperesrosen}).

Proposition \ref{propositionwrz2} implies therefore that a random walk on $\Z^2 \wr F$ is  individually almost invariant on the scale $f(n)$, for any $f(n)$ such that $(\ln f(n))^2 /\ln (n) \to 0$ as $n \to \infty$. In other words, 
for any $\epsilon_n$, tending to $0$ as $ n \to \infty$,  the simple random random walk on $\Z^2 \wr F$ is individually almost invariant on the scale
$ \exp(\epsilon_n  \sqrt{\ln(n)})$.
Therefore, 
the  radius of individual almost invariance of this random walk is greater than $ \exp( \sqrt{\ln(n)})$. In particular, this radius is strictly larger asymptotically than $\ln(n)$.

The Corollary follows therefore from the fact that $L(n) \sim n/\ln(n)$ for any simple random walk on $\Z^2 \wr F$ (see \cite{erschleruspehi} for a switch walk switch random walk, and a similar argument works for any simple random walk on this group).

\section{Transition probabilities of piecewise automatic groups} \label{piecewise}

The goal of this subsection is to prove Proposition  \ref{propositionbadinvariance}

Given a function $F(n)$,
we need to construct a group $G$ such that  that  for all $\epsilon, K$ the distributions $\mu^{n_i}$ of simple  random walks  on $G$ are not $(\epsilon, K)$ almost invariant with respect to balls of radius $F(n)$, for none of the choice of the subsequence $n_i$.

We will use  the construction from \cite{erschlerpiecewise}.

In this paper, we define and study "piecewise automatic groups", and we show (see Proposition 1 in \cite{erschlerpiecewise} and the proof of this proposition) that for  $\tau_1$ being the standard finite state automaton  for the first Gigorchuk group (possibly extended to some larger alphabet than $0$ and $1$) and $\tau_2$ being   a finite state automata, containing $e,a,b, c, d$ as its states, such that $a,b,c,d$ generate the free product $\mathbb{Z}/2\mathbb{Z} * (\mathbb{Z}/2\mathbb{Z}   + \mathbb{Z}/2\mathbb{Z})$, then the corresponding
"piecewise automatic groups with returns" has the following properties.

For any $t_1<T_1 <t_2 <T_2<t_3 ...$ one constructs a group $G(t_1, t_2, ... T_1, T_2, ...)$ with the following properties. There exist "comparison groups" $A(t_1, T_1, t_2, T_2, \dots t_i)$ and $B(t_1, T_1, t_2, T_2, \dots t_i, T_i)$, such that
\begin{enumerate}

\item all groups  $A(t_1, T_1, t_2, T_2, \dots t_i)$ are commensurable with a group which imbeds as a subgroup in a finite direct power of the the first Grigorchuk group $G_1$,

\item all groups  $B(t_1, T_1, t_2, T_2, \dots t_i, T_i)$ are commensurable with a group which surjects to the group $G_2$, generated by the automaton $\tau_2$,

\item furthermore, for some $\Psi(x)$ tending to $\infty$ as $x$ tends to $\infty$
the balls of radius $\Psi(T_i)$ in 
  $G(t_1, t_2, ... , T_1, T_2, ...)$  and $A(t_1, T_1, t_2, T_2, \dots t_i)$ coincide,

\item the balls of radius $\Psi(r_{i+1})$ in 
  $G(t_1, t_2, ... ,T_1, T_2, ...)$  and $B(t_1, T_1, t_2, T_2, \dots t_i, T_i)$ coincide.

\end{enumerate}

In particular, these properties imply that the constructed groups are of sub-exponential growth, so far as the sequence 
$t_1<T_1 <t_2 <T_2<t_3 ...$  grows quickly enough. In particular, symmetric finitely supported random walks on these group have trivial Poisson boundary. By Proposition \ref{propositiontrivialboundary} we know that in this case the group admits a non-trivial radius of almost invariance.

Observe that the groups $B(t_1, T_1, t_2, T_2, \dots t_i, T_i)$ are non-amenable, for any quickly enough growing sequence $t_1<T_1 <t_2 <T_2<t_3 <...$.  This implies that any non-degenerate random walk on these groups has non-trivial Poisson boundary.
By  Proposition \ref{propositiontrivialboundary} we know that for any $K,\epsilon$ there exists $C$ (depending on  $B(t_1, T_1, t_2, T_2, \dots t_i, T_i)$) such that the random walk on  $B(t_1, T_1, t_2, T_2, \dots t_i, T_i)$ is not $(K,\epsilon)$ on the balls of radius $C$.

The remark \ref{remarkCayley} below implies therefore that for any $F(n)$, tending to $\infty$ there exists  a sequence 
$t_1<T_1 <t_2 <T_2<t_3< ...$ such that the corresponding piecewise automatic group satisfies the following property.
For any $K, \epsilon$ and any sufficient large $n$
any almost invariance radius $r_{\epsilon,K}(n) \le F(n)$, and this completes the proof of the Proposition
\ref{propositionbadinvariance}.

\begin{remark} \label{remarkCayley}Let $A$ and $B$ the groups generated by $S_A$ and $S_B$ such that their marked Cayley graphs are isometric in the balls of radius $R\ge 1$. Let $\mu_A$ and $\mu_B$ be probability measures supported on $S_A$ and $S_B$ such that $\pi(\mu_A)=\mu_B$, for isomorphism $\pi$ between the balls of Cayley graphs (e.g. $\mu_A$ is an equidistribution on $S_A$ and $\mu_B$ is an equidistribution on $S_B$).
If the distributions of $(A,\mu_A)$ are $(K,\epsilon)$ almost invariant for $g\in A$ such that $ l_{S_A} \le r(n)$ and all $n$, then for any
$n$ such that $n+r(n)\le R$ the distributions of $(B,\mu_B)$ are $(K,\epsilon)$ almost invariant for $g \in  B$ such that $l_{S_B} \le r(n)$.

\end{remark}
{\bf Proof.} Observe that for any n  such that $n+r(n)\le R$, any $h\in A$ in the support of $\mu_A^{*n}$ and 
any $g \in A$ such that $l_A(g)\le r(n)$, it holds $l_A(gh)\le R$.

\section{Non constant limit distributions for ratios of transition probabilities}

\subsection{Limit distributions for word lengths}

Given a finitely generated group $G$ and a finite set of generators $S$ of $G$, we consider the word length
$l_S$. Given a sequence $g_n \in G$ (or a sequence of subsets $V_n \subset G$ and a sequence $\mu_n$ of probability measures on $G$, we can ask what is the limit behavior of $l_S(xg_n)-l_S(x)$ (where $g_n$ is a fixed sequence or $g_n\in V_n$), considered as a function on the probability space $(G,\mu_n)$. As before, we are interested in this paper in the case where $\mu_n= \mu^{*n}$ for some probability measure $\mu$ on $G$.

\begin{lemma} \label{lemmalengthfreegroups}[Free groups, distributions of $l_S(xh)-l_S(x)$]
Let $S$ be the free generating set of $G=F_m$, and $\mu$ be a measure which is equidistributed on $S \cup S^{-1}$. Fix an element $h \in F_m$ .

There exists $\epsilon_n$, tending to $0$ as $n\to \infty$ such that for any  $h\in G$ and 
 any $k \le l_S(h)$  it holds
$|\mu^{*n} (g: l_S(gh) < l_S(g)+l(h) -2 k) -  1/v_m(k)| \le \epsilon_n$. 
\end{lemma}

{\bf Proof.}
Let $l_S(g)=L$ and $l_S(h)=K$. Write $g= s_{i_1} \dots s_{i_L}$, $s_j \in S$ and 
$h= t_{j_1} \dots t_{i_K}$, $t_j \in S$. If $t_{j_1} \ne s_{i_L}^{-1}$, then $l_S(gh)=K+L$.
If $t_{j_1} = s_{i_L}^{-1}$ but $t_{j_2} \ne s_{i_{L-1}}^{-1}$, then $l_S(gh)=K+L-2$. 
If $t_{j_1} = s_{i_L}^{-1}$, $t_{j_2} = s_{i_{L-1}}^{-1}$, but $t_{j_3} \ne s_{i_{L-2}}^{-1}$
 then $l_S(gh)=K+L-3$ etc

Observe that $\mu^{*n} (x)= \mu^{*n}(y)$ for any $x,y: l_S(x)= l_S(y)$.

Since the random walk $(G,\mu)$ is transient, we know 
 that for any $C>0$ it holds
$\mu^{*n}(B(e,C)) \to 0$, where $B(e,C)$ is the ball of radius $C$ in the word length $l_S$.

Fix an element $h$ of length $L$. 
Observe that if we chose an element $g$ at random from the sphere of radius $K>L$ in $(G,l_S)$, then with probability $1-1/v_m(1)$ the last letter $s_{i_L}^{-1}$ of $g$ is not equal to the inverse of the first letter $t_{j_1}^{-1}$ of $h$. In this case$ l_S(gh) = l_S(g)+l(h)$
With probability $1-1/v_m(2)$ at least one two among the two pairs of letters are not equal:  either 
$t_{j_1} \ne s_{i_L}^{-1}$ or t $t_{j_2} \ne s_{i_{L-1}}^{-1}$. In this case $l_S(gh) \ge l_S(g)+l(h)-2$.
For any $r\ge 1$ the following holds: at least $m$ among the letters of $g$ and the corresponding inversed letters of $h$ are not equal with proability $1-1/v_m(r)$, and in this case $l_S(gh) \ge l_S(g)+l(h)-2r$

To prove the statement of the Lemma it suffices therefore to consider  $\epsilon_n= \mu^{*n}(B(e,K))$.

\subsection{Non-constant limit distributions}

Though no group with non-trivial boundary can admit a non-trivial scale for almost invariance, it may happen that such groups admit a scale  of  weak almost constant contraction. That is, a scale
on which the multiplication by an element $g$ multiplies the transition probability by some constant $C_g$,  up to a multiplicative constant and with probability close to one. This is the case for simple random walk on a free group $F_m$.

Let $v_m(n)=(2m-1) (2m)^{n-1}$ be the spherical growth function of the free group. with respect to the free generating set.

\begin{example} \label{examplefree}
Let $S$ be the free generating set of $F_m$, $m\ge 2$, and $\mu$ be a measure which is equidistributed on $S \cup S^{-1}$.

Let $\nu_m$ be the probability measure on the set $0, 1, 2, \dots $ such that 
$\nu_m(-x) = 1/v_m(x) - 1/v_m(x+1) = 1/2m (1/(2m-1)^{x}- 1/(2m-1)^{x+1})$ for any integer $x \ge 0$.
Consider the function $f(x)=-2x$ on this probability space.
Take any sequence $g_n \in F_m$ such that $l(g_n) \to \infty$ and  $l(g_n)/\sqrt{n} \to 0$, where $l$ is some word metric on $F_m$.
Then the distributions of the function $\mu^{*n}(g)/\mu^{*n}(g g_n)$  on the probability space $\mu^{*n}$ tend to that of $f$ on $\nu_m$.

\end{example}

{\bf Proof.}
Consider the random walk $X_n$  on $\mathbb{\Z^+}$, reflected at zero, which goes with probability
$1/(2m)$ to the left and $(2m-1)/2m$) to the right for any $x \in \mathbb{\Z^+}$, $x\ne 0$. 
Observe that the simple random walk on $F_m$ visits an element of length $l$ with probability of $X_n=l$, and that for two elements $g$ and $h$ in the free group of the same word length with respect to $S$ it holds $\mu^{*n}(g) =\mu^{*n}(h)$ for all $n$.

In particular,
$$
\mu^{^*n}(g) = 1/v(n) P[X_n=l_S(g)].
$$

\begin{remark}
Let $X_n^{o} \in \mathbb{Z}^+$ be such that $(X_n^{o} - (m-1)/m n))/\sqrt{n} \to 0$ and such that
$(X_n^{o} - [(m-1)/m n]))$ is divided by two. Then
$P[X_n= X_n^{o}]/ P[X_n = [(m-1)/m n]] \to 1$. Otherwise,  $(X_n^{o} - [(m-1)/m n]))$ is not divided by two for each $n$, then $P[X_n= X_n^{o}]/ P[X_n = [(m-1)/m n]+1] \to 1$. 

\end{remark}

The statement of Example \ref{examplefree} follows therefore from the Lemma \ref{lemmalengthfreegroups}.


\end{document}